\documentstyle{article}
\begin{document}

\author{Kosta Do\v sen}
\title{A Note on the Set-Theoretic Representation of Arbitrary Lattices}
\date{}
\maketitle

\begin{abstract}
Every lattice is isomorphic to a lattice whose elements are sets of sets and
whose operations are intersection and the operation $\vee ^{*}$ defined by
$A\vee ^{*}B=A\cup B\cup \{Z:(\exists X\in A)(\exists Y\in B)X\cap Y\subseteq
Z\}$. This representation spells out precisely Birkhoff's and Frink's
representation of arbitrary lattices, which is related to Stone's
set-theoretic representation of distributive lattices. ({\it AMS Subject
Classification, 1991: }06B15)
\end{abstract}

\noindent As a generalization of his representation theory for Boolean
algebras, Stone has developed in [4] a representation theory for
distributive lattices. This representation theory has set-theoretic and
topological aspects. Set-theoretically, every distributive lattice $L$ is
isomorphic to a set lattice $L^{*}$, i.e. a lattice whose elements are sets
and whose operations are intersection and union. In Stone's representation,
the elements of $L^{*}$ are certain subsets of the set $F(L)$ of prime
filters of $L$. Topologically, $F(L)$ can be viewed as a $T_0$-space with
the elements of $L^{*}$ constituting a subbasis.

Following ideas of Priestley's [3], Urquhart has developed in [5] the
topological aspects of this representation theory to cover arbitrary bounded
lattices. However, Birkhoff and Frink had already in [1] (section 6) a
simple set-theoretic representation for arbitrary lattices, also inspired by
Stone, but different from Urquhart's representation.

In the Birkhoff-Frink representation, every lattice $L$ is isomorphic to a
lattice $L^{*}$ whose elements are sets of sets, whose meet operation is
intersection and whose join operation is a set-theoretic operation $\vee
^{*} $ unspecified by Birkhoff and Frink. The elements of $L^{*}$ are
certain subsets of a set $F(L)$, which may be either the set of all filters
of $L$, or the set of all principal filters of $L$, or any set of filters of 
$L$ that for every pair of distinct elements of $L$ has a filter containing
one element of the pair but not the other. Stone's set-theoretic
representation for distributive lattices may be viewed as a special case of
the Birkhoff-Frink representation: if for a distributive lattice $L$ we take 
$F(L) $ to be the set of all prime filters of $L$, then $\vee ^{*}$
collapses into set-theoretic union.

The aim of this note is to make precise some details of the Birkhoff-Frink
representation, which doesn't seem to be very well known. We shall
explicitly characterize the operation $\vee ^{*}$ when $F(L)$ is the set of
all filters of $L$, or of all principal filters of $L$. The interest of this
exercise is in applications that may be found in the models of
nondistributive nonclassical logics, where the semantic clause for
disjunction may be derived from the operation $\vee ^{*}$.

Let $L=\langle D,\wedge ,\vee \rangle $ be an arbitrary lattice, and let
$F(L)=\{X:X$ is a filter of $L\}$. For every $a\in D$, let $f(a)=\{X\in
F(L):a\in X\}$. Let now $D^{*}=\{f(a):a\in D\}$, and let 
\begin{eqnarray*}
f(a)\wedge ^{*}f(b) &=&f(a)\cap f(b), \\
f(a)\vee ^{*}f(b) &=&f(a)\cup f(b)\cup \{Z\in F(L):(\exists X\in
f(a))(\exists Y\in f(b))X\cap Y\subseteq Z\}.
\end{eqnarray*}
The second of these equalities corresponds to the semantic clause for
disjunction introduced in [2] (section 3.2), which has since found its way
into a number of papers on models of substructural logics.

In the proof of the following proposition we assume for $a\in D$ that
$[a)=\{b\in D:a\leq b\}$; that is, $[a)$ is the principal filter generated by
$a$. \vspace{0.2cm}

\noindent {\bf Proposition 1.} {\it The following equalities hold}:

\vspace{0.2cm}

(1.1) $f(a)\wedge ^{*}f(b)=f(a\wedge b)$,

(1.2) $f(a)\vee ^{*}f(b)=f(a\vee b)$.

\vspace{0.2cm}

\noindent {\it Proof.} The proof of (1.1) is quite straightforward, and
we only need to consider the proof of (1.2). So suppose $Z\in f(a)\vee
^{*}f(b)$. If $a\in Z$ or $b\in Z$, then, since $Z$ is a filter, $a\vee b\in
Z$. If, on the other hand, for some $X$ and $Y$ we have that $a\in X$, $b\in
Y$ and $X\cap Y\subseteq Z$, then, since $X$ and $Y$ are filters, $a\vee
b\in X\cap Y $ , and so $a\vee b\in Z$. For the converse, suppose $Z\in
f(a\vee b)$, that is $a\vee b\in Z$. If $c\in [a\vee b)$, then $a\vee b\leq c
$, and, since $Z$ is a filter, $c\in Z$. So $[a\vee b)\subseteq Z$, but,
since $[a)\cap [b)=[a\vee b)$, we have that $[a)\cap [b)$ $\subseteq Z$.
Hence for some $X$, namely $[a)$, and some $Y$, namely $[b)$, we have that
$a\in X$, $b\in Y$ and $X\cap Y\subseteq Z$, and so we have proved (1.2).

\vspace{0.2cm}

Since it is quite easy to see that $f:D\rightarrow D^{*}$ is one-one and
onto, we obtain that $L=\langle D,\wedge ,\vee \rangle $ is isomorphic to
$L^{*}=\langle D^{*},\wedge ^{*},\vee ^{*}\rangle .$

Note that we obtain the isomorphism of $L$ with $L^{*}$ also when $F(L)$ is
taken to be the set of all principal filters of $L$, and not the set of all
filters of $L$. Another alternative, yielding again the isomorphism of $L$
with $L^{*}$, is to replace $\vee ^{*}$ by the operation $\vee ^{**}$
defined by 
\[
A\vee ^{**}B=\{Z:(\exists X\in A)(\exists Y\in B)X\cap Y\subseteq Z\}.
\]

We have preferred to work with $\vee ^{*}$, rather than with the more simply
defined operation $\vee ^{**}$, which coincides with $\vee ^{*}$ on $D^{*}$
as it was defined up to now, in order to be able to connect smoothly the
isomorphism of $L$ and $L^{*}$ with Stone's representation theory. This
connection is made by the following proposition.

\vspace{0.2cm}

\noindent {\bf Proposition 2.} {\it If} $L$ {\it is a distributive
lattice and} $F(L)$ {\it is the set of all prime filters of} $L$,
{\it then} $f(a)\vee ^{*}f(b)=f(a)\cup f(b)$.

\vspace{0.2cm}

\noindent {\it Proof. }Suppose $Z\in f(a)\vee ^{*}f(b)$. As in the
proof of the previous proposition, it follows that $a\vee b\in Z$. Since $Z$
is prime, $a\in Z$ or $b\in Z$, that is $Z\in f(a)\cup f(b)$. The converse,
namely, $f(a)\cup f(b)\subseteq $ $f(a)\vee ^{*}f(b)$, is trivial.

\vspace{0.2cm}

This trivial converse can, however, be blocked if $\vee ^{*}$ is replaced by 
$\vee ^{**}$. Indeed, suppose $a\in Z$; then we must show that for some
prime filters $X$ and $Y$ we have that $a\in X$, $b\in Y$ and $X\cap
Y\subseteq Z$. The prime filter $X$ can be $Z$, but, since $b$ may be the
least element of $L$, there is no guarantee that there is a prime, i.e.
proper, filter $Y$ such that $b\in Y$.

To conclude, we note that for the sake of symmetry we can define $f(a)\wedge
^{*}f(b)$ either as $f(a)\cap f(b)\cap \{Z\in F(L):(\exists X\in
f(a))(\exists Y\in f(b))X\cup Y\subseteq Z\}$, or as $\{Z\in F(L):(\exists
X\in f(a))(\exists Y\in f(b))X\cup Y\subseteq Z\}$; both of these sets are
equal to $f(a)\cap f(b)$. In these new definitions of $\wedge ^{*}$, unions
of filters occur where in the definitions of $\vee ^{*}$ and $\vee ^{**}$ we
had intersections. Then remark that the set of filters $F(L)$, which is a
semilattice with $\cap $, is not necessarily closed under $\cup $.

\vspace{0.5cm}

\begin{center}
{\bf References}
\end{center}

\vspace{0.2cm}

[1] G. Birkhoff and O. Frink, Jr, Representation of lattices by sets, 
{\it Trans. Amer. Math. Soc.} 64(1948), 299-316.

[2] K. Do\v sen, Sequent systems and groupoid models II, {\it Studia
Logica} 48(1989), 41-65.

[3] H.\thinspace A. Priestley,\ Representation of distributive lattices by
means of ordered Stone spaces, {\it Bull. London Math. Soc.} 2(1970),
186-190.

[4] M.\thinspace H. Stone, Topological representation of distributive
lattices and Brouwerian logics, {\it \v Casopis P\v est. Mat.} 67(1937),
1-25.

[5] A. Urquhart, A topological representation theory for lattices, {\it 
Algebra Universalis} 8(1978), 45-58.

\vspace{0.5cm}

{\it Matemati\v cki institut}

{\it Knez Mihailova 35, p.f. 367}

{\it 11000 Beograd, Jugoslavija}

{\it email: kosta@mi.sanu.ac.yu} \

\end{document}